\magnification=\magstep1
\nopagenumbers

\centerline{\bf Speed of convergence of two-dimensional Fourier  
integrals}
\vskip 2em
\centerline{Mark A Pinsky , Northwestern University}

\vskip 2em

\noindent {\bf 1. Introduction}

Recently [2,3] we found necessary and sufficient conditions for the  
convergence at a preassigned point of the spherical partial sums of  
the Fourier integral in a class of piecewise smooth functions in 
Euclidean space. These yield elementary  examples of divergent  
Fourier integrals in three dimensions and higher. Meanwhile, several  
years ago Gottlieb and Orsag[1] observed that in two dimensions we  
may expect slower convergence at certain points, specifically for 
Fourier-Bessel series of radial functions. 

In this paper we investigate the rate of convergence of the spherical  
partial sums of the Fourier integral for a class of piecewise smooth  
functions. The basic result is an asymptotic expansion which allows  
us to read off the rate of convergence at a pre-assigned point.

\vskip 2em
\noindent{\bf 2. Statement of results}
\vskip 1em

${\bf R}^2$ denotes the Euclidean plane and  $f\in L^1({\bf R}^2)$ is  
an integrable function with Fourier transform and spherical partial  
sum denoted by
$${\hat f}(\xi_1,\xi_2) =  
{1\over (2\pi)^2}\int_{{\bf  
R}^2}f(x_1,x_2)e^{-i(\xi_1x_1+\xi_2x_2)}dx_1dx_2,$$
$$f_M(x_1,x_2) = \int_{\xi_1^2 + \xi_2^2\le M^2}{\hat  
f}(\xi_1,\xi_2)e^{i(\xi_1x_1+\xi_2x_2)}d\xi_1 d\xi_2.$$
The spherical mean value with respect to $x=(x_1,x_2)$ is defined by
$${\bar f}_x(r) = {1\over 2\pi}\int_0^{2\pi} f(x_1 + r\cos \theta,  
x_2 + r\sin\theta)\,d\theta.$$
We say that $f$ is piecewise smooth  of class $C^k$ with respect to  
$x\in{\bf R}^2$ iff the mapping $r\rightarrow {\bar f}_x(r)$ is  
piecewise 

   $C^k$ and vanishes for $r\ge r(x)$.
 Such a function has left and right limits at all points, together  
with the requisite derivatives. These will agree except at a finite  
set of points which we denote by $a_i, 1\le i \le K.$ The jumps are  
denoted by

$$\delta {\bar f}_x(a_i) := {\bar f}_x (a_i + 0) - {\bar  
f}_x(a_i-0),$$
with a corresponding notation for the jumps of the derivatives.
\vskip 1em
\noindent {\bf Theorem} {\sl Suppose that} $f$ {\sl is piecewise 
smooth of class} $C^2$ {\sl with respect to} $x\in{\bf R}^2$. {\sl  
Then we have the asymptotic} ($M\rightarrow\infty$) {\sl formula}

$$\eqalign{f_M(x)-{\bar f}_x(0+0)&= \sqrt{2\over M\pi}\sum_{1\le i  
\le K}
 \delta {\bar f}_x(a_i) {\cos(Ma_i - \pi/4)\over \sqrt{a_i}}\cr
&+{1\over M}{\bar f'}_x(0+0)  +o(1/M).\cr}$$

\vskip 1em

\noindent{\bf Corollary} {\sl If the spherical mean value with  
respect to $x$ is continuous, then we have}
$$f_M(x) - {\bar f}_x(0+0) = O(1/M),\quad M\rightarrow \infty.$$
 Otherwise 

$$-\infty <\liminf_M {\sqrt M}[f_M(x) - {\bar f}_x(0+0) ]<
\limsup_M {\sqrt M}[f_M(x) - {\bar f}_x(0+0) ]<\infty.$$
\vskip 2em
\noindent {\bf Example} If we choose $f(x) =1_{[0,R]}(|x|)$ as the  
indicator function of a disc of radius $R$ centered at $(0,0)$,
 then $\delta {\bar f}_x(r) = 0$ unless $x=0,r=R$, where $\delta 
{\bar f}_0(R) =-1$; thus 

$$\eqalign{f_M(0,0) &=1-  \sqrt{2\over M\pi R}\cos(MR -\pi/4) + O(1/M),\cr
f_M(x_1,x_2) &=1+O(1/M)\qquad 0<x_1^2 + x_2^2 <R^2,\cr
f_M(x_1,x_2) &={1\over 2}+O(1/M)\qquad x_1^2 + x_2^2 =R^2,\cr
f_M(x_1,x_2) &= O(1/M)\qquad\qquad x_1^2 + x_2^2 >R^2.\cr}$$
The speed of convergence at the center is strictly slower than at all  
other points.
\vskip 2em
\noindent{\bf 3. Proofs}

The spherical partial sum is written directly in terms of the  
spherical mean value by writing [3]

$$f_M(x) = \int_0^a MJ_1(Mr) {\bar f}_x(r), \, dr\qquad a=r(x).$$
From the identities for Bessel functions, we have

$$MJ_1(Mr)=-{d\over dr} J_0(Mr),$$
so that we may integrate by parts:

$$\eqalign{f_M(x) &= -\int_0^a  {\bar f}_x(r) {d\over dr} J_0(Mr) 
dr\cr
&={\bar f}_x(0+0) + \sum_{1\le i\le K}\delta {\bar f}_x(a_i)J_0(Ma_i)  
+ 
\int_0^aJ_0(Mr){\bar f}'_x(r)\, dr.\cr}$$
The Bessel function has the asymptotic behavior  $J_0(x) =  
\sqrt{(2/\pi x)}[ \cos(x-\pi/4) +O(1/x)]$,

\noindent  $x\rightarrow\infty$, which yields the first term. To 
handle the integral term, we define $K(x) = \int_0^x J_0(t) \,dt$.  
From the asymptotic behavior of $J_0$ it follows that $K(x)$ is  
bounded with a limit $K(\infty) = \int_0^\infty J_0(t) \,dt=1$ (see  
the proof in the appendix below). Therefore we may integrate by parts  
once again to obtain

$$\eqalign{\int_0^aJ_0(Mr){\bar f'}_x(r)\, dr&= {1\over M}\int_0^a  
{\bar f'}_x(r) {d\over dr}
K(Mr)\, dr\cr
&=-{1\over M} \bigl(\sum_{1\le i \le K} \delta{\bar  
f'}_x(a_i)K(Ma_i)+ \int_0^a{\bar f''}_x(r) K(Mr)\, dr\bigr).\cr}$$
When $M\rightarrow \infty$ the first set of terms within the  
parenthesis clearly converges to

\noindent $K(\infty)\sum_{1\le i \le K} \delta{\bar f'}_x(a_i)$; the  
second term converges, by the dominated convergence theorem, to  
$K(\infty)\int_0^a{\bar f''}_x(r) \, dr$. Applying the fundamental  
theorem of calculus on each sub-interval $(a_{i-1},a_i)$ and  
collecting terms gives the stated form.
The proof is complete.

\vskip 2em

\noindent{\bf 4. Comparison with one-dimensional convergence}
\vskip 1em
In order to put the above in some perspective, we outline here a 
treatment of the corresponding questions for Fourier integrals in one  
variable, where we find 

a generic rate of convergence of $1/M$, even in the absence of  
continuity.

An integrable function $f(x), x\in {\bf R^1}$ has the Fourier  
transform

\noindent${\hat f}(\xi) = (2\pi)^{-1}\int_{{\bf R}^1} f(x) e^{-i\xi  
x} \,dx$ with spherical partial sum $f_M(x) = \int_{-M}^M {\hat  
f}(\xi) e^{i\xi x}\, d\xi.$
This can be expressed in terms of the {\it symmetric average} 

${\bar f}_x(t) = {1\over 2}[f(x+t) + f(x-t)]$ as

$$f_M(x) = {2\over \pi}\int_0^\infty {\sin Mt\over t}{\bar f}_x(t)\,  
dt$$ 

Using this representation, we now state and prove

\noindent {\bf Proposition.} Suppose that $t\rightarrow {\bar  
f}_x(t)$ is piecewise $C^2$ and zero for $t>a$. Then we have the 
asymptotic formula

$$f_M(x) -{{\bar f}_x}(0+0) = {2\over M \pi}[ {\bar f}_x'(0+0) +
\sum_{1\le i \le K}\delta {\bar f}_x(a_i) {\cos Ma_i\over a_i}]  
+o(1/M).$$

\noindent {\bf Proof.} We write
$$f_M(x) - {\bar f}_x(0+0) ={2\over \pi}\int_0^a {{\bar f}_x(t) - 
{\bar f}_x(0+0)\over t}\sin Mt\, dt
- {2\over \pi}{\bar f}_x(0+0)\int_a^\infty {\sin Mt\over t}\, dt,$$
where the last integral is conditionally convergent. In the first 
integral we write $g(t) = {{\bar f}_x(t) - {\bar f}_x(0+0)\over t}$  
and integrate by parts in the form
$$ \eqalign{\int_0^a g(t) \sin Mt\, dt &=-(1/M)\int_0^a g(t) d[\cos  
Mt]\cr
&=(1/M)[g(0) + \sum_{1\le i \le K}  \delta g(a_i) \cos Ma_i + 
\int_0^a g'(t) \cos Mt\, dt].\cr}$$
But $\delta g(a_i) = {\delta {\bar f}_x(a_i)\over a_i}$. 

The second integral can also be integrated by parts to obtain the 
asymptotic form
$$\int_a^\infty {\sin Mt\over t}\, dt = {\cos Ma\over Ma}+  
O(1/M^2).$$
Combining these produces the stated result.

\noindent {\bf Remark. } A parallel result hold for Fourier series on  
the circle $-\pi  <x<\pi$. In this case the spherical partial sum of  
$2M=2n+1$ terms is written

$$f_M(x) = {2\over \pi}\int_0^\pi {\sin Mt\over 2\sin(t/2)}{\bar 
f}_x(t)\, dt.$$
The above computations can be replicated with the result

$$f_M(x) -{{\bar f}_x}(0+0) = {2\over M \pi}[ {\bar f}_x'(0+0) +
\sum_{1\le i \le K}\delta {\bar f}_x(a_i) {\cos Ma_i\over  
2\sin(a_i/2)}] +o(1/M)$$

\noindent {\bf Conclusion.} We note that the one-dimensional result  
(for both Fourier series and integrals) differs from the  
two-dimensional result in the appearance of the derivative term  
${\bar f}_x'(0+0)$
at the {\it very first level} of the asymptotic analysis.

\vskip 2em

\noindent{\bf 5. Appendix. Computation of $K(\infty)$}
\vskip 1em
The Bessel function can be written
$$J_0(x)={1\over \pi}\int_{-1}^1 {\cos xt \over  
\sqrt{(1-t^2)}}\,dt.$$
Thus by the Fubini theorem 
$$\int_0^M J_0(x)\,dx = {1\over \pi}\int_{-1}^1 {\sin Mt\over t} 
{1\over \sqrt{(1-t^2)}}\,dt.$$
This is the partial Fourier integral of the integrable function  
$(1-t^2)^{-1/2}$ at $t=0$. 

When $M\rightarrow\infty$ we can appeal to the convergence of  
one-dimensional
Fourier integrals to conclude  
$$\lim_{M\rightarrow\infty}\int_0^M J_0(x) \,dx = 1.$$

\vfill\eject
\noindent{\bf References}
\vskip 2em
1. D. Gottlieb and S. Orsag, {\it Numerical Analysis of Spectral 
Methods}, SIAM, 1977.

\vskip 1em
2. M. Pinsky, Fourier inversion for piecewise smooth functions in 
several variables, Proceedings of the American Mathematical Society,  
118(1993), 903-910.
\vskip 1em

3. M. Pinsky, Pointwise Fourier inversion and related eigenfunction  
expansions, Communications of Pure and Applied Mathematics, 47(1994),  
653-681.

\bye